\documentclass[12pt]{article}
 \usepackage{latexsym}
 \usepackage{amssymb}
 \usepackage{graphicx}
 \usepackage{amsmath}

 \usepackage{dsfont}
 \usepackage{amssymb}
 \usepackage{booktabs}
 \usepackage{multirow}

\newtheorem{Theorem}{Theorem}

\newtheorem{Definition}{Definition}
\newtheorem{Proposition}{Proposition}

\newtheorem{Lemma}{Lemma}

\newtheorem{Algorithm}{Algorithm}

\newcommand{\A}{{\cal A}}
\newcommand{\B}{{\cal B}}

\newcommand{\I}{{\cal I}}

\newcommand{\x}{{\bf x}}
\newcommand{\y}{{\bf y}}
\newcommand{\z}{{\bf z}}
\newcommand{\cc}{{\bf c}}
\newcommand{\e}{{\bf e}}
\newcommand{\0}{{\bf 0}}

\newcommand{\w}{{\bf w}}

\newcommand{\qed}{\nobreak \ifvmode \relax \else
      \ifdim\lastskip<1.5em \hskip-\lastskip
      \hskip1.5em plus0em minus0.5em \fi \nobreak
      \vrule height0.75em width0.5em depth0.25em\fi}

\def \ep{\hbox{ }\hfill$\Box$}

\def \d{{\bf d}}

\begin{document}
\title{A Semismooth Newton Method for Tensor Eigenvalue Complementarity Problem}

\author{Zhongming Chen
\thanks{ School of Mathematical
Sciences and LPMC, Nankai University, Tianjin 300071, P.R. China.
Email: czm183015@gmail.com.
This author's work was done
when he was visiting The Hong Kong Polytechnic University.}, \quad
Liqun Qi
\thanks{ Department of Applied
Mathematics, The Hong Kong Polytechnic University, Hung Hom,
Kowloon, Hong Kong. Email: maqilq@polyu.edu.hk. This author's work was supported by the Hong
Kong Research Grant Council (Grant No. PolyU 502111, 501212, 501913 and 15302114).}
}

\date{\today}
\maketitle
\begin{abstract}
\noindent
In this paper, we consider the tensor eigenvalue complementarity problem which is closely related to the optimality conditions for polynomial optimization, as well as a class of differential inclusions with nonconvex processes. By introducing an NCP-function, we reformulate the tensor eigenvalue complementarity problem as a system of nonlinear equations. We show that this function is strongly semismooth but not differentiable, in which case the classical smoothing methods cannot apply. Furthermore, we propose a damped semismooth Newton method for tensor eigenvalue complementarity problem. A new procedure to evaluate an element of the generalized Jocobian is given, which turns out to be an element of the B-subdifferential under mild assumptions. As a result, the convergence of the damped semismooth Newton method is guaranteed by existing results. The numerical experiments also show that our method is efficient and promising.

\vspace{5mm}

\noindent {\bf Key words:}\hspace{2mm} Tensor eigenvalue complementarity problem, NCP-function, semismooth Newton method, B-subdifferential

\vspace{2mm}

\noindent {\bf AMS subject classifications (2010):}\hspace{2mm}
15A18; 15A69; 90C33; 49M15


\end{abstract}

\newpage
\section{Introduction}
\hspace{14pt}
Complementarity problems has developed into a very fruitful discipline in the field of mathematical programming, which were originally studied in the standard optimality conditions for linear and smooth nonlinear optimization. The distinguishing feature of a complementarity problem is a set of complementarity conditions. Each of these conditions requires that the product of two or more quantities should be zero. They appear prominently in the study of equilibria problems and arise naturally in numerous applications from economics, engineering and the sciences \cite{FP}.

As a special type of complementarity problems, the matrix eigenvalue complementarity problem has also been well-studied for the last decade, which first appeared in the stability analysis of finite dimensional mechanical systems with frictional contact \cite{DFJM}.
Mathematically speaking, given a matrix $A \in \mathbb{R}^{n\times n}$ and a positive definite matrix $B \in \mathbb{R}^{n\times n}$, the matrix eigenvalue complementarity problem \cite{JSR,QJH} consists of finding a scalar $\lambda > 0$ and a vector $\x \in \mathbb{R}^n \setminus \{\0\}$ such that
$$
\x \geq \0 , \quad (\lambda B - A) \x \geq \0,  \quad  \big\langle \x, (\lambda B - A) \x \big\rangle =0 ,
$$
where $\langle \cdot, \cdot \rangle$ denotes the inner product of two vectors.
For more details about matrix eigenvalue complementarity problems and their applications, the readers are referred to \cite{AR,DMFJ,DS,JRRS,JSRR}.

Meanwhile, as a natural extension of the concept of matrices, a real $m$th order $n$-dimensional tensor $\A=(a_{i_1 \ldots i_m})$ is a multidimensional array with each entry $a_{i_1 \ldots i_m} \in \mathbb{R}$ for any $ i_1, \ldots, i_m \in [n]$, where $[n]=\{ 1,2, \ldots, n\}$. Denote the set of all real $m$th order $n$-dimensional tensors by $T_{m, n}$. For any vector $\x \in \mathbb{R}^n$,
let $\A \x^{m-1}$ be a vector in $\mathbb{R}^n$ whose $i$th component is defined by
$$
(\A \x^{m-1})_i =\sum_{i_2, \ldots, i_m=1}^n a_{i i_2 \ldots i_m} x_{i_2} \ldots x_{i_m}.
$$
And let $\A \x^m$ be the scale denoted by $\A \x^m= \x^\top \A \x^{m-1}$,
which is exactly a homogeneous polynomial of $\x$ with degree $m$.
We say a tensor $\A$ is {\it symmetric} if its entries are invariant under permutation.
Denote the set of all real symmetric $m$th order $n$-dimensional tensors by $S_{m, n}$.
A tensor $\A \in S_{m,n}$ is called {\it positive definite} if $\A \x^m > 0$ for all $\x \neq 0$.
Clearly, when $m$ is odd, there is no positive definite tensors.

It is possible that the ideas of eigenvalues of tensors had been raised earlier.
However, it was the independent work of Lim \cite{Lim} and Qi \cite{Qi} that initiated the rapid developments of the spectral
theory of tensors. Moreover, these definitions can all be unified under generalized tensor eigenpair framework as follows, introduced by Chang, Pearson and Zhang \cite{CPZ}. Let $\A, \B \in T_{m,n}$. Assume further that $m$ is even and $\B$ is positive definite. We say $(\lambda, \x) \in \mathbb{C} \times \{ \mathbb{C}^n \setminus \{ \0 \} \}$ is a generalized eigenpair if
$$
\A \x^{m-1} = \lambda \B \x^{m-1}.
$$
Different choices of $\B$ yield different versions of the tensor eigenvalue problem \cite{CPZ,KM}.
After that, the study of tensors and the spectra of tensors with their various applications has attracted extensive attention and interest.

In this paper, we consider the tensor eigenvalue complementarity problem (TEiCP), which consists of finding a scalar $\lambda >0$ and a vector
$\x \in \mathbb{R}^n \setminus \{ \0 \}$ such that
\begin{equation}\label{e1}
\x \geq \0 , \quad (\lambda \B - \A) \x^{m-1} \geq \0,  \quad  \big\langle \x, (\lambda \B - \A) \x^{m-1} \big\rangle =0 ,
\end{equation}
where $\A \in T_{m,n}$ and $\B \in S_{m,n}$ is positive definite.
This problem is not only a natural generalization of the classical eigenvalue complementarity problem for matrices,
but also closely related to the optimality conditions for polynomial optimization \cite{SQ2}, a kind of nonlinear differential dynamical system \cite{CYY},
as well as a class of differential inclusions with nonconvex processes \cite{LHQ1}.
In particular, without the constraint $\lambda >0$,
any solution $\lambda$ of (\ref{e1}) is called a Pareto H-eigenvalue or a Pareto Z-eigenvalue by Song and Qi \cite{SQ2}
when tensor $\B$ has different special forms, respectively. The properties of Pareto eigenvalues and their relationship
with polynomial optimization are also studied in \cite{SQ2}. The general properties of TEiCP, including the solution
existence and uniqueness, are systematically investigated by Chen, Yang and Ye \cite{CYY}. Note that the eigenvalue
complementarity problem under a given index set is also considered in \cite{CYY}.
When the nonnegative cones in (\ref{e1}) are replaced
by a closed convex cone and its dual cone, TEiCP is called the cone eigenvalue
complementarity problem for high-order tensors by Ling, He and Qi \cite{LHQ1}. Moreover, as a natural extension of
quadratic eigenvalue complementarity problem for matrices, Ling, He and Qi \cite{LHQ2} also consider the high-degree
eigenvalue complementarity problem for tensors. Some properties of Pareto eigenvalues are further studied in \cite{XL}.
Another kind of complementarity problems related to tensors is considered in \cite{CQW, SLQX, LQX, SQ1}.

On the other hand, some algorithms for computing the solutions of TEiCP have been proposed, such as
shifted projected power method \cite{CYY}, scaling-and-projection algorithm \cite{LHQ1} and alternating direction
method of multipliers \cite{LHQ2}. Notice that all these methods are first order algorithms that are based
on gradient information. In this paper, we present a semismooth Newton method for computing the solutions of TEiCP.
It turns out that TEiCP is a parameterized nonlinear complementarity problem. By introducing an NCP-function to
get rid of the nonnegative constraints in (\ref{e1}), we reformulate it as a system of nonsmooth operator equations.
The main difficulty is that this function is not differentiable. As a result, the classical Newton method cannot apply.
Fortunately, this function is strongly semismooth and the semismooth Newton method has been well studied since the work
of Qi and Sun \cite{QS}. In order to implement the semismooth Newton method, we also propose a new procedure to
evaluate an element of the generalized Jacobian, which turns out to be an element of the B-subdifferential under some
mild assumptions. The numerical results indicate that our method is efficient to compute the solutions of TEiCP.

The rest of this paper is organized as follows. In Section 2, we recall some basic definitions and results in
nonsmooth analysis and nonlinear complementarity problem. In Section 3, we reformulate TEiCP as a system of nonlinear
equations. In Section 4, we present a damped semismooth Newton method for TEiCP, as well as convergence analysis.
Some numerical results are reported in Section 5. In Section 6, we give two final remarks.

Throughout this paper, we assume that $m$ is even and $\B \in S_{m,n}$ is positive definite.
We use small letters $x, y, \ldots, $ for scalers, small bold letters $\x, \y, \ldots, $ for vectors,
capital letters $A, B, \ldots, $ for matrices, calligraphic letters $\A, \B, \ldots, $ for tensors.
All the tensors discussed in this paper are real.


\section{Preliminaries}
\hspace{14pt}
In this section, we present some basic definitions and properties in nonsmooth analysis and nonlinear complementarity problem,
which will be used in the sequel.

Suppose that $F: U \subseteq \mathbb{R}^{n_1} \rightarrow \mathbb{R}^{n_2}$ is a locally Lipschitz function,
where $U$ is nonempty and open. By Rademacher's
Theorem, $F$ is differentiable almost everywhere. Let $D_F \subseteq \mathbb{R}^{n_1}$ denote the set of points at which $F$ is differentiable.
For any $\x \in D_F$, we write $JF(\x)$ for the usual $n_2 \times n_1$ Jacobian matrix of partial derivatives.
The {\it B-subdifferential} of $F$ at $\x \in U$ is the set defined by
$$
\partial_B F(\x) := \left\{  V \in \mathbb{R}^{n_2 \times n_1} :  \exists \{\x^k \} \subseteq D_F \text{ with } \x^k \rightarrow \x, JF(\x^k) \rightarrow V \right\}
$$
And Clark's {\it generalized Jacobian} of $F$ at $\x$ is the set defined by
$$
\partial F(\x) = \text{co}(\partial_B F(\x)),
$$
where ``co'' denotes the convex hull. In the case $n_2 =1$, $\partial F(\x)$ is called the {\it generalized gradient}.
Some fundamental properties about generalized Jacobian are given as follows.
For more details, one can refer to \cite{Cla}.

\begin{Proposition}\label{p1}
Suppose that the function $F: U \subseteq \mathbb{R}^{n_1} \rightarrow \mathbb{R}^{n_2}$ is locally Lipschitz,
where $U$ is nonempty and open. Then for any $\x \in U$, we have
\begin{itemize}
\item[(a)] $\partial F(\x)$ is a nonempty convex compact subset of $\mathbb{R}^{n_2 \times n_1}$;
\item[(b)] $\partial F(\x) = \partial_B F(\x) = \{ JF(\x)\}$ if $F$ is continuously differentiable at $\x$;
\item[(c)] $\partial F(\x) \subseteq \partial f^1(\x) \times \partial f^2 (\x) \times \cdots \times \partial f^m (\x)$,
where $F(\x)= [ f^1(\x), f^2(\x), \ldots, f^m(\x)]$ and the latter denotes the set of all matrices
whose $i$th row belongs to $\partial f^i(\x)$ for each $i$.
\end{itemize}
\end{Proposition}

Let $U \subseteq \mathbb{R}^{n_1}$ be nonempty and open.
The function $F : U \rightarrow \mathbb{R}^{n_2} $ is {\it semisoomth} \cite{QS} at $\x \in \mathbb{R}^{n_1}$, if
it is locally Lipschitz at $\x$ and if
$$
\lim_{\substack{V \in \partial F(\x+t \tilde{\d}) \\ \tilde{\d} \rightarrow \d, \ t\downarrow 0}} V \tilde{\d}
$$
exists for all $\d \in \mathbb{R}^n$. If $F$ is semismooth at all $\x \in U$, we call $F$ semismooth on $U$.
And the function $F$ is called {\it strongly semismooth} \cite{QY} if it is semismooth and for any $\x \in U$ and $V \in \partial F(\x+\d )$,
$$  V \d - F'(\x; \d) = O( \| \d \|^2 ) , \quad  \d \rightarrow \0 ,$$
where $F'(\x; \d)$ denotes the directional derivative \cite{BCS} of $F$ at $\x$ in direction $\d$, i.e.,
$$
F'(\x; \d) = \lim_{t \downarrow 0} { F(\x+ t \d) - F(\x)  \over t }.
$$
It is worth mentioning \cite{QS} that if the function $F$ is semismooth, the directional derivative $F'(\x; \d)$ exists for all $\d \in \mathbb{R}^n$
and
$$ F'(\x; \d) =  \lim_{\substack{V \in \partial F(\x+t \tilde{\d}) \\ \tilde{\d} \rightarrow \d, \ t\downarrow 0}} V \tilde{\d} .$$

Suppose that $F : \mathbb{R}^n \rightarrow \mathbb{R}^n$ is directional differentiable at $\x$, we say that $F$ is {\it BD-regular} \cite{Qi}
at $\x$ if for any $\d \in \mathbb{R}^n \setminus \{ \0\}$,
$$ F'(\x ; \d) \neq \0 .$$
We say that $F$ is {\it strongly BD-regular} at $\x$ if all $V \in \partial_B F(\x)$ are nonsingular.
These concepts can be very important for the convergence analysis of semismooth Newton method.

In the following, the classical nonlinear complementarity problem is introduced. And it will be shown in the next section that
tensor eigenvalue complementarity problem is a special kind of parameterized nonlinear complementarity problem.

\begin{Definition}
Given a mapping $F: \mathbb{R}^n_+ \rightarrow \mathbb{R}^n $, the nonlinear complementarity problem, denoted by NCP($F$),  is to find a vector
$\x \in \mathbb{R}^n$ satisfying
$$ \x \geq \0, \quad  F(\x) \geq \0, \quad  \langle x, F(\x) \rangle  = 0  .$$
\end{Definition}

Many solution methods developed for NCP or related problems are based on reformulating
them as a system of equations using so-called NCP-functions. Here, a function
$\phi : \mathbb{R}^2 \rightarrow \mathbb{R}$ is called an {\it NCP-function} if
$$ \phi(a, b)=0 \Leftrightarrow  ab = 0, \ a \ge  0, \ b \ge 0. $$
Given an NCP-function $\phi$, let us define
\begin{equation}\label{e2}
 \Phi (\x) = \big[ \phi(x_i, F_i(\x) )\big]_{i=1}^n  .
\end{equation}
By definition, $\x \in \mathbb{R}^n$
is a solution of NCP($F$) if and only if it solves the system of equations $\Phi (\x) = \0$.

Here, we present some NCP-functions which are widely used in nonlinear complementarity problems.
For more details about NCP-functions and their smoothing approximations, one can refer to \cite{QSZ,SQ,YLQ,ZCT} and references therein.
\begin{itemize}
\item The min function \cite{SS}
$$ \phi_{min}(a,b) := a - (a-b)_+ .$$
\item The Fischer-Burmeister function \cite{Fis}
$$ \phi_{FB}(a,b):= (a+b) - \sqrt{a^2+b^2} .$$
\item The penalized Fischer-Burmeister function \cite{CCK}
$$ \phi_{\tau}(a,b):= \tau \phi_{FB}(a,b) + (1-\tau) a_+ b_+ ,$$
where $\tau \in (0,1)$.
\end{itemize}
Here, $x_+ =\max\{x,0\}$ for $x \in \mathbb{R}$. It has been shown that all these NCP-functions are
globally Lipschitz continuous, directionally differentiable, and strongly semismooth.
And their generalized gradients are given as follows.

\begin{Proposition}\label{p2}
Let $\phi_{min}(a,b)$, $\phi_{FB}(a,b)$ and $\phi_{\tau}(a,b)$ de defined as above. Then
\begin{itemize}
\item[(a)] The generalized gradient $\partial \phi_{min}(a,b)$ is equal to the set of all $(v_a, v_b)$ such that
$$ (v_a,v_b) = \left\{ \begin{array}{ll}
(1,0)   & \text{if } a < b, \\
(1-v,v) & \text{if } a = b, \\
(0,1)   & \text{if } a > b,
\end{array}
\right.  $$
where $v$ is any scalar in the interval $[0,1]$.
\item[(b)] The generalized gradient $\partial \phi_{FB}(a,b)$ is equal to the set of all $(v_a, v_b)$ such that
$$ (v_a,v_b) = \left\{ \begin{array}{ll}
\big(1-\frac{a}{\|(a,b)\|},1-\frac{b}{\|(a,b)\|} \big)    & \text{if } (a,b) \neq (0,0), \\
(1-\sigma,1-\eta) & \text{if } (a,b) = (0,0),
\end{array}
\right.  $$
where $(\sigma,\eta)$ is any vector satisfying $\|(\sigma, \eta) \| \leq 1$.
\item[(c)] For any $\tau \in (0,1)$, the generalized gradient $\partial \phi_{\tau}(a,b)$ is equal to the set of all $(v_a, v_b)$ such that
$$ (v_a,v_b) = \left\{ \begin{array}{ll}
\tau \big(1-\frac{a}{\|(a,b)\|},1-\frac{b}{\|(a,b)\|} \big) + (1-\tau) (b_+ \partial a_+, a_+ \partial b_+ )   & \text{if } (a,b) \neq (0,0), \\
\tau(1-\sigma,1-\eta) & \text{if } (a,b) = (0,0),
\end{array}
\right.  $$
where $(\sigma,\eta)$ is any vector satisfying $\|(\sigma, \eta) \| \leq 1$ and
$$ \partial x_+ = \left\{ \begin{array}{ll}
0     & \text{if } x < 0, \\
{[0,1]} & \text{if } x = 0, \\
1     & \text{if } x > 0.
\end{array}
\right.  $$
\end{itemize}
\end{Proposition}

\section{Reformulation}
\hspace{14pt}
Suppose that $m$ is even. As mentioned before, the tensor eigenvalue complementarity problem has the form of
\begin{equation}\label{e3}
\text{(TEiCP): Find } \lambda >0, \x \neq 0 \text{ such that } \left\{ \begin{array}{l}  \w=(\lambda \B -\A)\x^{m-1}  \\ \w \geq 0 \\  \x \geq 0  \\ \w^\top \x =0 , \end{array} \right.
\end{equation}
where $\A \in T_{m,n}$ and $\B \in S_{m,n}$ is positive definite.
Note that any solution with $\w =0$ is a generalized tensor eigenpair of $(\A, \B)$.
Denote the solution set of (\ref{e1}) by $\sigma(\A, \B)$, i.e.,
$$ \sigma(\A,\B) = \big\{ (\lambda, \x) \in \mathbb{R}_{++} \times \mathbb{R}^n\setminus \{ \0 \} :  \0 \leq \x \perp (\lambda \B - \A) \x^{m-1} \geq \0 \big\}.$$
Notice that if $(\lambda, \x) \in \sigma (\A,\B)$, then $(\lambda, s \x) \in \sigma (\A, \B)$ for any $s >0$.
Without loss of generality, we only consider the solution satisfying $\|\x\|_2 =1$.
On the other hand, it is clear that
$$\sigma(\A, \B) = \sigma(\tilde{\A},\B),$$
where $\tilde{\A} \in T_{m,n}$ is the unique semi-symmetric tensor \cite{NQ} such that $\A \x^{m-1} = \tilde{\A} \x^{m-1}$ for all $\x \in \mathbb{R}^n$.
Hence, we always assume that $\A \in T_{m,n}$ is semi-symmetric.

By introducing a new variable $t \in \mathbb{R}$, we denote
\begin{equation}\label{e4}
F(\x,t):= (t^2 \B - \A ) \x^{m-1}, \quad \ \forall \ \x \in \mathbb{R}^n, t \in \mathbb{R}.
\end{equation}
As a result, TEiCP can be regarded as a parameterized nonlinear complementarity problem,
i.e.,
\begin{equation}\label{e5}
\x \ge \0, \quad  F(\x,t) \ge \0, \quad  \x^\top F(\x,t) = 0,
\end{equation}
with the constraint $\x^\top \x = 1 $.
It follows that the TEiCP can be represented compactly by
the system of nonlinear equations
\begin{equation}\label{e6}
H(\z) = \0 ,
\end{equation}
where $\z =(\x ,t)\in \mathbb{R}^{n+1}$, $t \neq 0$, and $H :   \mathbb{R}^{n+1} \rightarrow \mathbb{R}^{n+1}$ is defined by
\begin{equation}\label{e7}
H(\z) = \left(   \begin{array}{c} \Phi(\z)  \\ \x^\top \x -1   \end{array}  \right) ,
\end{equation}
where the mapping $\Phi: \mathbb{R}^{n+1} \rightarrow \mathbb{R}^n$ is given by
\begin{equation}\label{e8}
 \Phi (\z) = \big[ \phi(x_i, F_i(\x,t) )\big]_{i=1}^n  ,
\end{equation}
and $\phi(a,b)$ is an NCP-function. Moreover, a natural metric function of $H(\z)$ is given by
\begin{equation}\label{e9}
\Psi (\z) = \frac{1}{2} H(\z)^\top H(\z).
\end{equation}

By using the techniques in nonlinear complementarity problem, we can see that finding a solution of TEiCP is equivalent to solve
the corresponding system of nonlinear equations.
\begin{Proposition}\label{p3}
Let $H(\z)$ be defined by (\ref{e7}). If $(\lambda, \x)$ is a solution of (\ref{e3}) with $\| \x \|=1 $, then
$H(\z) =0$ with $\z=(\x, \pm \sqrt{\lambda}) \in \mathbb{R}^{n+1}$. On the other hand, if $H(\z) = 0$ with
$z=(\x, t) \in \mathbb{R}^{n+1}$ and $t \neq 0$, then $(t^2, \x)$ is a solution of (\ref{e3}) with $\| \x \|=1 $.
\end{Proposition}


Let $H_{min}(\z)$, $H_{FB}(\z)$ and $H_{\tau}(\z)$ be the functions defined by (\ref{e7}), corresponding to the NCP-functions
$\phi_{min}$, $\phi_{FB}$ and $\phi_{\tau}$, respectively. In the following, we will show that these functions are all strongly semismooth,
which can be very important for nonsmooth Newton methods.
\begin{Lemma}\label{l1}
The functions $\Phi_{min}(\z)$, $\Phi_{FB}(\z)$ and $\Phi_{\tau}(\z)$ are strongly semismooth,
where $\Phi_{min}(\z)$, $\Phi_{FB}(\z)$ and $\Phi_{\tau}(\z)$ are defined by (\ref{e8}), corresponding to the NCP-functions
$\phi_{min}$, $\phi_{FB}$ and $\phi_{\tau}$, respectively.
\end{Lemma}
{\bf Proof.} Notice that for any $\z =(\x,t)\in \mathbb{R}^{n+1}$, the function $F(\z)= (t^2 \B - \A ) \x^{m-1}$ is continuously differentiable and its Jacobian $JF(\z)$ is locally Lipschitz continuous. It follows from Theorem 1 of \cite{SQ} that $\Phi_{FB}(\z)$ and $\Phi_{\lambda}(\z)$ are strongly semismooth. Similarly, $\Phi_{min}(\z)$ is strongly semismooth since the composition of strongly semismooth functions is again strongly semismooth \cite{Mif}.
\ep

\begin{Theorem}\label{t1}
The functions $H_{min}(\z)$, $H_{FB}(\z)$ and $H_{\tau}(\z)$ are
strongly semismooth. Moreover, for any $\z \in \mathbb{R}^{n+1}$, we have
\begin{equation}\label{e10}
\partial H(\z) \subseteq \left\{  \left( \begin{array}{cc} \multicolumn{2}{c}{V}   \\  2\x^\top & 0 \end{array} \right) \in \mathbb{R}^{(n+1)\times(n+1)} : V \in \partial\Phi(\z)   \right\}.
\end{equation}
\end{Theorem}
{\bf Proof.} It is clear that $H_{n+1} (\z)= \x^\top \x -1 $ is continuously differentiable. By Lemma \ref{l1},
the functions $\Phi_{min}(\z)$, $\Phi_{FB}(\z)$ and $\Phi_{\tau}(\z)$ are strongly semismooth. It follows that
$H_{min}(\z)$, $H_{FB}(\z)$ and $H_{\tau}(\z)$ are strongly semismooth since all their components are strongly semismooth.
Moreover, by Proposition 2.6.2 of \cite{Cla}, (\ref{e10}) holds immediately.
\ep

\section{Semismooth Newton Method}
\hspace{14pt}
In order to establish a semismooth Newton method for TEiCP, we need to obtain en element of $\partial H(\z)$.
First, we have the following result.
\begin{Proposition}\label{p4}
Suppose that the mapping $F: \mathbb{R}^{n_1 + n_2} \rightarrow \mathbb{R}^{n_1}$ is continuously differentiable and
the function $\phi : \mathbb{R}^2 \rightarrow \mathbb{R}$ is locally Lipschitz.
Let $ \Phi : \mathbb{R}^{n_1 + n_2} \rightarrow \mathbb{R}^{n_1}$ be the mapping such that
$$ \Phi (\z) = \left[ \phi(x_i, F_i(\z))\right]_{i=1}^{n_1} , \quad  \forall \ \z =(\x , \y) \in \mathbb{R}^{n_1+ n_2} .$$
Then for any $\z \in \mathbb{R}^{n_1 +n_2}$, we have
$$
\partial \Phi (\z) \subseteq \left( D_a(\z), \0_{n_1 \times n_2} \right) +  D_b(\z) JF(\z),
$$
where $D_a(\z) = diag\{ a_i (\z) \} $ and $D_b(\z) = diag \{ b_i(\z) \}$ are diagonal matrices in $\mathbb{R}^{n_1 \times n_1}$
with entries $(a_i (\z), b_i(\z)) \in \partial \phi (x_i, F_i(\z))$, where $\partial \phi (x_i, F_i(\z))$ denotes the set
$\partial \phi(a,b)$ with $(a,b)$ being replaced by $(x_i, F_i(\z))$.
\end{Proposition}

Let $F(\z)$ be defined in (\ref{e4}) and $\phi (a,b)$ be one of the NCP-functions $\phi_{FB}$ and $\phi_{\tau}$.
Since $\A \in T_{m,n}$ is semi-symmetric, by simple computation, the Jacobian of $F$ at $\z$ is given by
$$
JF(\z) = \left[ (m-1) (t^2 \B -\A) \x^{m-2} \quad  2 t \B \x^{m-1} \right] \in \mathbb{R}^{n \times (n+1)}.
$$
By Propositions \ref{p2} and \ref{p4}, we can obtain the overestimation of $\partial \Phi_{FB} (\z)$ and $\partial \Phi_{\tau} (\z)$, respectively.
In the following, we present a procedure to obtain an element of $\partial \Phi_{\tau} (\z)$
for any $\z \in \mathbb{R}^{n+1}$, where $\tau \in (0, 1]$. Note that $\partial \Phi_{1} (\z) = \partial \Phi_{FB} (\z)$.

\begin{Algorithm}\label{alg1}
{\rm A procedure to generate an element $V \in \partial\Phi_{\tau}(\z)$}

\noindent {\bf Step 0.} Given $\tau \in (0,1]$, $\z =(\x,t) \in \mathbb{R}^{n+1} $ and let $V_i$ be the $i$th row of a matrix $V \in \mathbb{R}^{n\times (n+1)}$.

\noindent {\bf Step 1.} Set $S_1 = \{ i \in [n] : x_i =0,  F_i(\x,t) =0\}$, $S_2 = \{ i \in [n] : x_i =0, F_i(\x,t) >0\}$, $S_3 = \{ i \in [n] : x_i >0,  F_i(\x,t) =0\}$ and $S_4 = \{ i \in [n] : x_i >0, F_i(\x,t) >0\}$.

\noindent {\bf Step 2.} Let $\cc \in \mathbb{R}^n $ such that $c_i =1$ for $i \in S_1 \cup S_2 \cup S_3$ and $0$ otherwise.

\noindent {\bf Step 3.} For $i \in S_1$, set
$$
V_i = \tau\left(1 + \frac{c_i}{\| (c_i, \nabla_{\x} F_i(\z)^\top \cc) \|} \right) (\e_i^\top , 0) + \tau\left(1 +\frac{\nabla_{\x} F_i(\z)^\top \cc}{\| (c_i,  \nabla_{\x} F_i(\z)^\top \cc) \|} \right) \nabla F_i(\z)^\top.
$$

\noindent {\bf Step 4.} For $i \in S_3$, set
$$
V_i = \left\{ \begin{array}{ll}
\big(\tau +(1-\tau)x_i \big) \nabla F_i(\z)^\top   & \text{ if } \  \nabla_{\x} F_i(\z)^\top \cc < 0  , \vspace{2mm} \\
\tau \nabla F_i(\z)^\top                  & \text{ otherwise.}
\end{array} \right.
$$

\noindent {\bf Step 5.} For $i \in S_4$, set
\begin{multline*}
V_i = \left[ \tau\left(1- \frac{x_i}{\| (x_i, F_i(\z) )\|}\right) +(1-\tau)F_i(\z)  \right] (\e_i^\top, 0)  \\
+\left[ \tau\left(1- \frac{F_i(\z)}{\| (x_i, F_i(\z) )\|}\right) +(1-\tau)x_i \right] \nabla F_i (\z)^\top.
\end{multline*}

\noindent {\bf Step 6.} For $i \not \in S_1 \cup S_3 \cup S_4$, set
$$
V_i = \tau \left(1- \frac{x_i}{\| (x_i, F_i(\z) )\|}\right) (\e_i^\top, 0) + \tau\left(1- \frac{F_i(\z)}{\| (x_i, F_i(\z) )\|}\right) \nabla F_i (\z)^\top.
$$
\end{Algorithm}

\medskip
For $\tau \in (0,1)$, based on the overestimate of $\partial \Phi_{\tau}(\z)$, we can see that under some assumptions, the matrix
$$ G = \left( \begin{array}{cc} \multicolumn{2}{c}{V}   \\  2\x^\top & 0 \end{array} \right)$$
is an element in the B-subdifferential of $H_\tau $ at $\z$, where
$V \in \mathbb{R}^{n \times (n+1)}$ is the matrix generated by Algorithm \ref{alg1} with $\tau \in (0,1)$.

\begin{Theorem}\label{t2}
Let $\z =(\x,t) \in \mathbb{R}^{n+1}$ be given and let $V \in \mathbb{R}^{n \times (n+1)}$ be the matrix generated by Algorithm \ref{alg1}
with $\tau \in (0,1)$.
Suppose that $\nabla_{\x} F_i (\z)^\top \cc \neq  0$ for all $i \in S_3$. Then the matrix
$G = \left( \begin{array}{cc} \multicolumn{2}{c}{V}   \\  2\x^\top & 0 \end{array} \right)$
is an element of $\partial_B H_{\tau} (\z)$.
\end{Theorem}
{\bf Proof.} Notice that $H_{\tau}(\z)$ is differentiable except the set
$$\Omega: = \big\{ \z =(\x,t) \in \mathbb{R}^{n+1} : x_i \geq 0, F_i(\z) \geq 0, x_i F_i(\z) =0 \text{ for some } i \in [n] \big\}. $$
We shall generate a sequence $\{ \z^k \}_{k=1}^{\infty} \subseteq \mathbb{R}^{n+1}\setminus \Omega $ such that $JH(\z^k)$ tends to
the matrix $G$. Then the conclusion follows immediately by the definition of B-subdifferential.

The conclusion is trivial if $\z \not \in \Omega$, i.e., $S_1 \cup S_2 \cup S_3 = \emptyset$.
In the following, we suppose that $\z \in \Omega$, i.e., $S_1 \cup S_2 \cup S_3 \neq \emptyset$.
Let $\z^k = \z - \frac{1}{k} (\cc^\top, 0)$, where $\cc \in \mathbb{R}^n$ is the vector given in Step 2.
It is clear that $z^k_i < 0$ for $i \in S_1 \cup S_2$.
For $i \in S_1 \cup S_3$, by Taylor-expansion, we have
\begin{equation}\label{e11}
F_i(\z^k) = F_i(\z) + \nabla F_i (\zeta^k )^\top (\z^k - \z) = -\frac{1}{k} \nabla_{\x} F_i (\zeta^k)^\top \cc ,
\end{equation}
where $\zeta^k \rightarrow \z$ as $k \rightarrow \infty$.
Since $\nabla_{\x} F_i (\z)^\top \cc \neq  0$ for all $i \in S_3$, by continuity,
we have that for all $i \in S_3$, $F_i(\z^k) \neq 0$ when $k$ is large enough.
Hence, there exists $N >0$ such that $H_{\tau}(\z^k)$ is differentiable for all $k > N$.

For $i \in [n+1]$, let $ JH(\z^k)_i$ be the $i$th row of $JH(\z^k)$.
If $i \not \in S_1 \cup S_2 \cup S_3$ or $i = n+1$, by continuity, it is obvious that $JH(\z^k)_i$ tends to the $i$th row
of $G$. For $i \in S_1 \cup S_2$, we have
$$
JH(\z^k)_i = \tau \left(1- \frac{x_i^k}{\| (x_i^k, F_i(\z^k) )\|}\right) (\e_i^\top, 0) + \tau\left(1- \frac{F_i(\z^k)}{\| (x_i^k, F_i(\z^k) )\|}\right) \nabla F_i (\z^k)^\top.
$$
For $i \in S_3$, it is not difficult to show that
$$
JH(\z^k)_i = \left\{ \begin{array}{ll}
{\small \begin{split}
\left[ \tau\left(1- \frac{x_i^k}{\| (x_i^k, F_i(\z^k) )\|}\right) +(1-\tau)F_i(\z^k)  \right] (\e_i^\top, 0)  \\
& \hspace{-70mm}  +\left[ \tau\left(1- \frac{F_i(\z^k)}{\| (x_i^k, F_i(\z^k) )\|}\right) +(1-\tau)x_i^k \right] \nabla F_i (\z^k)^\top
\end{split}}                & \text{ if } \  \nabla_{\x} F_i(\z^k)^\top \cc < 0  , \vspace{4mm} \\
\tau \left(1- \frac{x_i^k}{\| (x_i^k, F_i(\z^k) )\|}\right) (\e_i^\top, 0) + \tau\left(1- \frac{F_i(\z^k)}{\| (x_i^k, F_i(\z^k) )\|}\right) \nabla F_i (\z^k)^\top            & \text{ if } \  \nabla_{\x} F_i(\z^k)^\top \cc > 0 .
\end{array} \right.
$$
Note that for $i \in S_1$, by substituting (\ref{e11}), we have
$$
\lim_{k \rightarrow \infty} \frac{x_i^k}{\| (x_i^k, F_i(\z^k) )\|} = \lim_{k \rightarrow \infty} \frac{-1/k}{\sqrt{(1/k)^2 + F_i(\z^k)^2}} =\frac{-1}{\sqrt{1 + (\nabla_{\x} F_i(\z)^\top \cc)^2}}.
$$
Similarly,
$$
\lim_{k \rightarrow \infty} \frac{F_i(\z^k)}{\| (x_i^k, F_i(\z^k) )\|} =\frac{-\nabla_{\x} F_i(\z)^\top \cc}{\sqrt{1 + (\nabla_{\x} F_i(\z)^\top \cc)^2}}.
$$
It follows that for $i \in S_1 \cup S_2 \cup S_3$, $JH(\z^k)_i$ tends to the $i$th row of the matrix $G$.
\ep

\medskip
We also mention that $H_{FB}(\z)$ is differentiable except the set
$$
\{ \z=(\x,t) \in \mathbb{R}^{n+1}: x_i=0, F_i(\z) =0 \text{ for some } i \in [n] \} .
$$
Hence, $H_{FB}(\z)_i$ is differentiable for all $i \in S_3$. By a similar proof of Theorem \ref{t2},
we can see that the matrix $G = \left( \begin{array}{cc} \multicolumn{2}{c}{V}   \\  2\x^\top & 0 \end{array} \right)$
is exactly an element of $\partial_B H_{FB} (\z)$ without any assumption, where $V \in \mathbb{R}^{n \times (n+1)}$ is the matrix generated by Algorithm \ref{alg1} with $\tau =1$.

\begin{Theorem}\label{t3}
Let $\z =(\x,t) \in \mathbb{R}^{n+1}$ be given and let $V \in \mathbb{R}^{n \times (n+1)}$ be the matrix generated by Algorithm \ref{alg1}
with $\tau =1$.
Then the matrix
$G = \left( \begin{array}{cc} \multicolumn{2}{c}{V}   \\  2\x^\top & 0 \end{array} \right)$
is an element of $\partial_B H_{FB} (\z)$.
\end{Theorem}

\medskip
Now we present some properties of the metric functions $\Psi_{FB}(\z)$ and $\Psi_{\tau}(\z)$.
\begin{Theorem}\label{t4}
The metric function $\Psi_{FB}(\z)$ and $\Psi_{\tau}(\z)$ defined in (\ref{e9}) is continuously differentiable with $ \nabla \Psi (\z) = G^\top H(\z)$ for any $G \in \partial H(\z)$.
\end{Theorem}
{\bf Proof.}
By known rules on the calculus of generalized gradients (see \cite{Cla}, Theorem
2.6.6), it holds that $\partial \Psi(\z) = H(\z)^\top \partial H(\z)$. Since it is easy to check that $H(\z)^\top \partial H(\z)$ is
single valued everywhere because the zero components of $H(\z)$ cancel the ``multivalued rows'' of $\partial H(\z)$ , we have by the corollary to Theorem 2.2.4 in \cite{Cla} that $\Psi(\z)$ is
continuously differentiable.
\ep

\medskip
Note that the metric function $\Psi_{min}(\z)$ is not continuously differentiable, which makes
the generalized Newton direction not necessarily a descent direction. Now we present a damped Newton method
for tensor eigenvalue complementarity problem. Here,
we only take the NCP-functions $\phi_{FB}$ and $\phi_{\tau}$.

\begin{Algorithm}\label{alg2}
{\rm Damped semismooth Newton method for TEiCP}

\noindent {\bf Step 0.} Given $\epsilon >0$, $\rho >0$, $p>2$, $\beta \in (0, \frac{1}{2})$ and choose $\z^0 =(\x^0, t^0)\in \mathbb{R}^{n+1}$. Set $k=0$.

\noindent {\bf Step 1.} If $\| H(\z^k) \| \leq \epsilon $, stop. Otherwise, go to Step 2.

\noindent {\bf Step 2.}
Compute
$$
G_k = \left( \begin{array}{cc} \multicolumn{2}{c}{V_k}   \\ 2(\x^k)^\top & 0 \end{array} \right),
$$
where $V_k \in \mathbb{R}^{n \times (n+1)} $ is the element of $\partial \Phi (\z^k)$ generated by Algorithm \ref{alg1}.
Find the solution $\d^k$ of the system
\begin{equation}\label{e12}
G_k \d = - H(\z^k).
\end{equation}
If $G_k$ in (\ref{e12}) is ill-conditioned or if the condition
$$
\nabla \Psi(\z^k)^\top \d^k \leq -\rho \| \d^k\|^p
$$
is not satisfied, set $\d^k = - \nabla \Psi (\z^k)$.

\noindent {\bf Step 3.} Find the smallest $i_k= 0,1, \ldots $ such that $\alpha_k = 2^{-i_k}$ and
$$ \Psi (\z^k + \alpha_k \d^k) \leq \Psi (\z^k)  + \beta \alpha_k \nabla \Psi(\z^k)^\top \d^k .$$
Set $\z^{k+1} = \z^k + \alpha_k \d^k$.

\noindent {\bf Step 4.} Set $k=k+1$ and go back to Step 1.
\end{Algorithm}

The global convergence of Algorithm \ref{alg2} is also guaranteed by the following theorem,
whose proof can be found in the papers \cite{DFK,Pa,Qi}.

\begin{Theorem}\label{t5}
Suppose that the solution set $\sigma(\A,\B)$ is nonempty. Let $\{\z^k\} \subseteq \mathbb{R}^{n+1}$ be generated by Algorithm \ref{alg2}.
Assume that $H(\z^k) \neq \0$ for all $k$. Then the conclusions (a) and (b) hold:
\begin{itemize}
\item[(a)] $\| H(\z^{k+1}) \| \leq \| H(\z^k) \|$ ;
\item[(b)] each accumulation point $\z^*$ of the sequence $\{ \z^k \}$ is a stationary point of $\Psi$.
Furthermore, if $H(\z)$ is strongly BD-regular at $\z^*$, then $\z^*$ is a zero of $H(\z)$ if and only if $\{ \z^k\}$ converges to $\z^*$ quadratically
and $\alpha_k$ eventually becomes 1. On the other hand, $\z^*$ is not a zero of $H(\z)$ if and only if $\{ \z^k\}$ diverges or
$\lim_{k \rightarrow \infty}\alpha_k =0$.
\end{itemize}
\end{Theorem}

Here, we make several remarks. First, if $H(\z^k) = \0$ for some $k$, Algorithm \ref{alg2} terminates right then with a zero solution of $H(\z)$.
Second, by Theorem \ref{t1}, $H(\z)$ is strongly semismooth. It follows from Lemma 2.3 of \cite{Qi} that
the directional differential $H'(\cdot, \cdot)$ is semicontinuous of degree 2 at $\z^*$.
Then the convergence is quadratic if $H(\z)$ is strongly BD-regular at $\z^*$. Third, when $H(\z)$ is not strongly BD-regular at $\z^*$, i.e.,
there exists a singular matrix $V \in \partial_B H(\z^*)$, the local convergence is also established by using adaptive constructs of outer inverses
\cite{CNQ}.

\section{Numerical results}
\hspace{14pt}
In this section, we present the numerical performance of the damped semismooth Newton method for the tensor eigenvalue complementarity problem.
All codes were written by using Matlab Version R2012b and the Tensor Toolbox Version 2.5 \cite{BK}.
And the numerical experiments were done on a laptop with an Intel Core i5-2430M CPU (2.4GHz) and RAM of 5.58GB.

In the implementation of Algorithm \ref{alg2}, we set the parameters $\epsilon = 10^{-6}$, $\rho = 10^{-10} $, $p=2.1$ and $\beta = 10^{-4}$.
We choose the penalized Fischer-Burmeister function with $\tau =0.95$. Numerically, we say that the matrix $G_k$ in (\ref{e12}) is ill-conditioned if $\kappa(G_k) \geq 10^{10}$, where $\kappa(M)= {\sigma_{max}(M) \over \sigma_{min} (M)}$ denotes the condition number of the matrix $M$.
We let Algorithm 2 run until any of the following situations occur:
\begin{itemize}
\item[(a)] $k=1000$                 \ \ \ \ \  \  \qquad      Failure (converge to a stationary point but not a solution),
\item[(b)] $\| H(\z)\| \leq \epsilon$   \ \ \ \qquad          Success (a solution has been detected).
\end{itemize}
For simplicity, the positive definite tensor $\B \in S_{m,n}$ is chosen by the identity tensor \cite{KM}, i.e.,
$\B \x^{m-1} = \x$ for all $\x^\top \x =1$ where $m$ is even.

Our first numerical experiment concerns the symmetric tensor $\A \in S_{6,4}$ described in Table 1 of \cite{CYY}.
We compare the numerical performance between damped semismooth Newton method and shifted projected power method \cite{CYY}.
It has been shown that this problem is solvable, i.e., the solution set is nonempty, since $A$ is symmetric with $a_{111111} > 0$ \cite{CYY}.
In order to get all possible solutions, we take 1000 random initial points for both methods. We adopt the following strategy:
one generates a random vector $\x^0 \in \mathbb{R}^4$ with each element uniformly distributed on the open interval $(0,1)$, a scalar $t^0 \in \mathbb{R}$
drawn from the standard normal distribution $N(0,1)$, and then normalizes $\x^0$ such that $\| \x^0 \| =1$.
Note that for the shifted projected power method, the initial point $\x^0$ needs to be chosen such that $\A (\x^0)^m > 0$ while
the damped semismooth Newton method do not have this constraint. We also record the time of finding the proper initial point for the shifted projected power method.

The numerical results of these two methods are reported in Tables \ref{tb1} and \ref{tb2}, respectively.
In Tables \ref{tb1} and \ref{tb2}, {\it No.} denotes the number of each solution detected by the method within 1000 random initial points.
{\it Ite.} denotes the average number of iteration for each solution. In Table \ref{tb1}, {\it Time1} and {\it Time2} denote the average time
of finding the initial point and the average time of iteration in second, respectively.  In Table \ref{tb2}, {\it Time} denotes the average time of iteration
and the values in bold are the solutions detected only by the damped semismooth Newton method.

Our second numerical experiment consist of applying the damped semismooth Newton method to a sample of 10 randomly generated tensors.
Given order $m$ and dimension $n$, the tensor $\A \in S_{m,n}$ is generated as \cite{CYY}, i.e.,
we select random entries from $[-1, 1]$ and symmetrize the result. To make TEiCP solvable, we reset its first entry by 0.5.
The idea is measuring the rate of success of the damped semismooth Newton algorithm with a given number of initial points.
The damped semismooth Newton algorithm is declared successful if a solution is found while working with a prescribed number of initial points
(for instance, 1, 5, or 10). The outcome of this experiment is reported in Table \ref{tb3}.

The third numerical experiment focuses on nonnegative tensors. As pointed out in \cite{CYY}, there exists a unique solution
for TEiCP when $\A \in T_{m,n}$ is irreducible nonnegative and $\B = \I \in S_{m,n}$, where $\I$ is the diagonal tensor with diagonal entries 1
and 0 otherwise. To be specific, this unique solution is exactly the largest H-eigenvalue of $\A$, associated with the unique positive eigenvector.
We generate a symmetric nonnegative tensor $\A \in S_{6,4}$ with all entries randomly selected from $(0,1)$. The tensor $\A$ is given in Table \ref{tb4}.
We test the damped semismooth Newton method with the initial point $\x^0 = \e / \| \e \|$ and $t^0 = \sqrt{\A (\x^0)^6 / \I (\x^0)^6}$, where
$\e= (1, 1, \ldots, 1)^\top \in \mathbb{R}^n$. The iteration of Algorithm \ref{alg2} is given in Table \ref{tb5}.
Here we also report the value of $\alpha_k$ derived in Step 3 of Algorithm \ref{alg2}.

It is worth mentioning that the damped semismooth Newton method can also work for nonsymmetric tensors while the shifted projected power method
does not work. We test the damped semismooth Newton method for randomly generated nonnegative tensors with all entries selected from the interval $(0,1)$.
Note that these nonnegative tensors are not symmetric in general, and the order $m$ is not necessary to be even. We use the same initial point in the third numerical experiment. Interestingly, we find that the damped semismooth Newton method always converges to the unique solution. We summarize our numerical results in Table \ref{tb6}. For each case, we use a sample of 100 random tensors to record the number of success ({\it Suc.}), the average number of iteration ({\it Ite.}), the average time of iteration ({\it Time}) and the average value of $\lambda$-solution ($\lambda^*$).

\begin{table}
\tabcolsep 1mm
\caption{All possible solutions detected by the shifted projected power method for the tensor $\A$ given in Table 1 of \cite{CYY}} \label{tb1}
\begin{center}
\begin{tabular}{c c l l c c c}
\toprule
No.   & $\lambda^*$   & \multicolumn{1}{l}{$\x^*$}         & \multicolumn{1}{l}{$\w^*$}         & Time1(s) & Ite.  &  Time2(s)  \\ \hline
73    & 0.5081	      & $(0.8158,	0, 0, 0.5784)^\top$    & $(0, 0.2647, 0.2295, 0)^\top$      & 0.011182 & 9.22  & 0.074547   \\
11    & 0.6136	      & $(0, 0, 0, 1)^\top$                & $(0.3403, 0.1706, 0.0559, 0)^\top$ & 0.009655 & 2.36  & 0.021255   \\
196   & 0.8181	      & $(0, 0.7804, 0.6253, 0)^\top$      & $(0.0735, 0, 0, 0.0243)^\top$      & 0.012117 & 9.08  & 0.078052   \\
402   & 0.8568	      & $(0, 0.8251, 0.5146, 0.2333)^\top$ & $(0.3184, 0, 0, 0)^\top$           & 0.011359 & 40.4  & 0.338430   \\
318   & 1.1666	      & $(0.5781, 0, 0.816, 0)^\top$       & $(0, 0.3347, 0, 0.4207)^\top$      & 0.013218 & 48.6  & 0.409752   \\
\bottomrule
\end{tabular}
\end{center}
\end{table}

\begin{table}
\tabcolsep 1.05mm
\caption{All possible solutions detected by the damped semismooth Newton method for the tensor $\A$ given in Table 1 of \cite{CYY}} \label{tb2}
\begin{center}
\begin{tabular}{c c l l c c}
\toprule
No.   & $\lambda^*$   & \multicolumn{1}{l}{$\x^*$}               & \multicolumn{1}{l}{$\w^*$}          & Ite.  &  Time(s)  \\ \hline
161   & {\bf 0.2655}  & ${\bf (0.2922, 0.4238, 0.8573, 0)^\top}$ & $(0, 0, 0, 0.3808)^\top$            & 31.06 &	0.782480  \\
44    & {\bf 0.2985}  & ${\bf (0.9738, 0, 0, 0.2274)^\top}$      & $(0, 0.261, 0.1069, 0)^\top$        & 38.14 &	1.058434  \\
174   & 0.5081	      & $(0.8158,	0, 0, 0.5784)^\top$          & $(0, 0.2647, 0.2295, 0)^\top$       & 7.11  &  0.139725  \\
64    & 0.6136	      & $(0, 0, 0, 1)^\top$                      & $(0.3403, 0.1706, 0.0559, 0)^\top$  & 7.05  &	0.126869  \\
24    & {\bf 0.8161}  & ${\bf (0, 0.7884, 0.6146, 0.0274)^\top}$ & $(0.1047, 0, 0, 0)^\top$            & 6.83  &	0.119179  \\
49    & 0.8181	      & $(0, 0.7804, 0.6253, 0)^\top$            & $(0.0735, 0, 0, 0.0243)^\top$       & 6	 &  0.106406  \\
253   & 0.8568	      & $(0, 0.8251, 0.5146, 0.2333)^\top$       & $(0.3184, 0, 0, 0)^\top$            & 8.23  &	0.189993  \\
135   & 1.1666	      & $(0.5781, 0, 0.816, 0)^\top$             & $(0, 0.3347, 0, 0.4207)^\top$       & 7.19  &	0.141976  \\
96    &   \multicolumn{5}{c}{converge to a stationary point but not a solution}                                       \\
\bottomrule
\end{tabular}
\end{center}
\end{table}

\begin{table}
\tabcolsep 5mm
\caption{The success rate of the damped semismooth Newton method} \label{tb3}
\begin{center}
\begin{tabular}{c c c c c}
\toprule
    &      & \multicolumn{3}{l}{Number of random initial points}  \\ \cline{3-5}
m   &   n  & 1        & 5           &   10          \\ \hline
4   &   5  &   70\%   &  100\%      &   100\%       \\
4   &  10  &   60\%   &  100\%      &   100\%       \\
4   &  20  &   30\%   &  90\%       &   100\%       \\
4   &  30  &   10\%   &  70\%       &   90\%        \\
4   &  40  &   10\%   &  50\%       &   90\%        \\
6   &   5  &   90\%   &  100\%      &   100\%       \\
6   &  10  &   80\%   &  100\%      &   100\%       \\
8   &   4  &   50\%   &  100\%      &   100\%       \\
8   &   5  &   50\%   &  100\%      &   100\%       \\

\bottomrule
\end{tabular}
\end{center}
\end{table}

\begin{table}
\caption{A random symmetric nonnegative tensor $\A =(a_{i_1 i_2 \ldots i_6})\in T_{6,4}$} \label{tb4}
\begin{center}
\scalebox{0.9}[0.9]{%
\begin{tabular}{l@{ = } c  l@{ = } c  l@{ = } r l@{ = } c}
\toprule
$a_{111111}$ & $  0.1197$, & $a_{111112}$ & $0.4859$, & $a_{111113}$ & $0.4236$, & $a_{111114}$ & $0.1775$,  \\
$a_{111122}$ & $  0.4639$, & $a_{111123}$ & $0.4951$, & $a_{111124}$ & $0.5322$, & $a_{111133}$ & $0.4219$,  \\
$a_{111134}$ & $  0.4606$, & $a_{111144}$ & $0.4646$, & $a_{111222}$ & $0.4969$, & $a_{111223}$ & $0.4649$,  \\
$a_{111224}$ & $  0.5312$, & $a_{111233}$ & $0.5253$, & $a_{111234}$ & $0.4635$, & $a_{111244}$ & $0.4978$,  \\
$a_{111333}$ & $  0.5562$, & $a_{111334}$ & $0.5183$, & $a_{111344}$ & $0.4450$, & $a_{111444}$ & $0.4754$,  \\
$a_{112222}$ & $  0.4992$, & $a_{112223}$ & $0.5420$, & $a_{112224}$ & $0.4924$, & $a_{112233}$ & $0.5090$,  \\
$a_{112234}$ & $  0.4844$, & $a_{112244}$ & $0.5513$, & $a_{112333}$ & $0.5040$, & $a_{112334}$ & $0.4611$,  \\
$a_{112344}$ & $  0.4937$, & $a_{112444}$ & $0.5355$, & $a_{113333}$ & $0.4982$, & $a_{113334}$ & $0.4985$,  \\
$a_{113344}$ & $  0.4756$, & $a_{113444}$ & $0.4265$, & $a_{114444}$ & $0.5217$, & $a_{122222}$ & $0.2944$,  \\
$a_{122223}$ & $  0.5123$, & $a_{122224}$ & $0.4794$, & $a_{122233}$ & $0.5046$, & $a_{122234}$ & $0.4557$,  \\
$a_{122244}$ & $  0.5332$, & $a_{122333}$ & $0.5161$, & $a_{122334}$ & $0.5236$, & $a_{122344}$ & $0.5435$,  \\
$a_{122444}$ & $  0.5576$, & $a_{123333}$ & $0.5685$, & $a_{123334}$ & $0.5077$, & $a_{123344}$ & $0.5138$,  \\
$a_{123444}$ & $  0.5402$, & $a_{124444}$ & $0.4774$, & $a_{133333}$ & $0.6778$, & $a_{133334}$ & $0.4831$,  \\
$a_{133344}$ & $  0.5030$, & $a_{133444}$ & $0.4865$, & $a_{134444}$ & $0.4761$, & $a_{144444}$ & $0.3676$,  \\
$a_{222222}$ & $  0.1375$, & $a_{222223}$ & $0.5707$, & $a_{222224}$ & $0.5440$, & $a_{222233}$ & $0.5135$,  \\
$a_{222234}$ & $  0.5770$, & $a_{222244}$ & $0.6087$, & $a_{222333}$ & $0.5075$, & $a_{222334}$ & $0.4935$,  \\
$a_{222344}$ & $  0.5687$, & $a_{222444}$ & $0.5046$, & $a_{223333}$ & $0.5226$, & $a_{223334}$ & $0.4652$,  \\
$a_{223344}$ & $  0.5289$, & $a_{223444}$ & $0.4810$, & $a_{224444}$ & $0.5310$, & $a_{233333}$ & $0.6187$,  \\
$a_{233334}$ & $  0.5811$, & $a_{233344}$ & $0.4811$, & $a_{233444}$ & $0.4883$, & $a_{234444}$ & $0.4911$,  \\
$a_{244444}$ & $  0.4452$, & $a_{333333}$ & $0.1076$, & $a_{333334}$ & $0.6543$, & $a_{333344}$ & $0.4257$,  \\
$a_{333444}$ & $  0.5786$, & $a_{334444}$ & $0.5956$, & $a_{344444}$ & $0.4503$, & $a_{444444}$ & $0.3840$.  \\
\bottomrule
\end{tabular}}
\end{center}
\end{table}

\begin{table}
\caption{The iteration of Algorithm \ref{alg2} for the tensor $\A$ given in Table \ref{tb4}} \label{tb5}
\begin{center}
\begin{tabular}{c c l l l }
\toprule
$k$ &$\lambda^k$ &$\x^k$                                  & $\| H(\z^k)\|$      &  $\alpha_k$  \\ \hline
1   & 514.2207   &$(0.4977, 0.5013, 0.5005, 0.5005)^\top$ &  1.30e-01           &    1         \\
2   & 515.3181   &$(0.4982, 0.5012, 0.5003, 0.5003)^\top$ &  1.05e-02           &    1         \\
3   & 515.4096   &$(0.4982, 0.5012, 0.5003, 0.5003)^\top$ &  1.08e-04           &    1         \\
4   & 515.4105   &$(0.4982, 0.5012, 0.5003, 0.5003)^\top$ &  1.20e-08           &    1         \\
5   & 515.4105   &$(0.4982, 0.5012, 0.5003, 0.5003)^\top$ &  1.83e-14           &    1         \\
\bottomrule
\end{tabular}
\end{center}
\end{table}

\begin{table}
\tabcolsep 5mm
\caption{Numerical results for random nonnegative tensors} \label{tb6}
\begin{center}
\begin{tabular}{cccccc}
\toprule
m   &   n  &  Suc.  & Ite.   &  Time(s)    & $\lambda^*$    \\  \hline
3   &  20  &  100   & 5.48   &  0.0998     &   0.2000e+03   \\
3   &  40  &  100   & 6.00   &  0.1375     &   0.8001e+03   \\
3   &  60  &  100   & 6.00   &  0.2215     &   1.8004e+03   \\
3   &  80  &  100   & 6.00   &  0.2904     &   3.2000e+03   \\
3   & 100  &  100   & 6.00   &  0.4424     &   4.9999e+03   \\
4   &  10  &  100   & 5.06   &  0.0851     &   0.4997e+03   \\
4   &  20  &  100   & 5.57   &  0.1542     &   4.0002e+03   \\
4   &  30  &  100   & 5.92   &  0.3636     &   1.3498e+04   \\
4   &  40  &  100   & 5.98   &  0.8690     &   3.1999e+04   \\
4   &  50  &  100   & 6.00   &  1.9368     &   6.2500e+04   \\
5   &   5  &  100   & 4.65   &  0.0677     &   0.3127e+03     \\
5   &  10  &  100   & 5.07   &  0.1137     &   5.0000e+03   \\
5   &  15  &  100   & 5.29   &  0.3102     &   2.5311e+04   \\
5   &  20  &  100   & 5.66   &  1.0001     &   7.9993e+04   \\
6   &   4  &  100   & 4.51   &  0.0695     &   0.5112e+03     \\
6   &   6  &  100   & 4.82   &  0.0968     &   3.8880e+03   \\
6   &   8  &  100   & 4.99   &  0.1695     &   1.6386e+04   \\
6   &  10  &  100   & 5.10   &  0.3610     &   5.0000e+04   \\
8   &   4  &  100   & 4.41   &  0.0873     &   8.1914e+03   \\
8   &   5  &  100   & 4.69   &  0.1826     &   3.9070e+04   \\
\bottomrule
\end{tabular}
\end{center}
\end{table}

\section{Final remarks}
\hspace{14pt}
In this paper, we propose a damped semismooth Newton method for tensor eigenvalue complementarity problem. Here we make two final remarks.

1. Given an index set $J \subseteq [n]$, the generalized tensor eigenvalue complementarity problem $\text{(TEiCP)}_J$ is also considered in \cite{CYY}.
In fact, we can also apply our damped semismooth Newton method to $\text{(TEiCP)}_J$. Since the results are similar, we omit them in this paper.

2. From the numerical results, we may give a new way to compute the largest H-eigenvalue of irreducible nonnegative tensors since this problem can be reformulated as a tensor eigenvalue complementarity problem equivalently.


\end{document}